\documentclass{amsart}

\theoremstyle{definition}

\theoremstyle{remark}

\reversemarginpar

\begin{document}

\title[A exp-log family]{The integrals in Gradshteyn and Ryzhik. Part 3: \\
Combinations of logarithms and exponentials.}

\author{Victor H. Moll}
\address{Department of Mathematics,
Tulane University, New Orleans, LA 70118}
\email{vhm@math.tulane.edu}

\subjclass{Primary 33}

\date{\today}

\keywords{Integrals}

\begin{abstract}
We present the evaluation of a family of exponential-logarithmic integrals.  
These have integrands of the form $P(e^{tx}, \ln x)$ 
where $P$ is a 
polynomial. The examples presented here appear in sections  $4.33, \, 4.34$
and $4.35$ in the classical 
table of integrals by I. Gradshteyn and I. Ryzhik.
\end{abstract}

\maketitle

\newcommand{\nn}{\nonumber}
\newcommand{\ba}{\begin{eqnarray}}
\newcommand{\ea}{\end{eqnarray}}
\newcommand{\ift}{\int_{0}^{\infty}}
\newcommand{\ione}{\int_{0}^{1}}
\newcommand{\ifft}{\int_{- \infty}^{\infty}}
\newcommand{\no}{\noindent}
\newcommand{\Ftwo}{{{_{2}F_{1}}}}
\newcommand{\realpart}{\mathop{\rm Re}\nolimits}
\newcommand{\imagpart}{\mathop{\rm Im}\nolimits}

\newtheorem{Definition}{\bf Definition}[section]
\newtheorem{Thm}[Definition]{\bf Theorem} 
\newtheorem{Example}[Definition]{\bf Example} 
\newtheorem{Lem}[Definition]{\bf Lemma} 
\newtheorem{Note}[Definition]{\bf Note} 
\newtheorem{Cor}[Definition]{\bf Corollary} 
\newtheorem{Prop}[Definition]{\bf Proposition} 
\newtheorem{Problem}[Definition]{\bf Problem} 
\numberwithin{equation}{section}

\section{Introduction} \label{intro} 
\setcounter{equation}{0}

This is the third in a series of papers dealing with the evaluation of 
definite integrals in the 
table of Gradshteyn and Ryzhik \cite{gr}. We consider here problems of the 
form
\begin{equation}
\ift e^{- tx} \, P(\ln x) \, dx, 
\label{class-0}
\end{equation}
\noindent
where $t >0$ is a parameter and $P$ is a polynomial. In future work we 
deal with the finite interval case
\begin{equation}
\int_{a}^{b} e^{- tx} \, P(\ln x) \, dx, 
\label{class-1}
\end{equation}
\noindent
where $a, \, b \in \mathbb{R}^{+}$ with $a < b$ and 
$t \in \mathbb{R}$. The classical example
\begin{equation}
\int_{0}^{\infty} e^{-x} \ln x \, dx = - \gamma, 
\label{class1}
\end{equation}
\noindent
where $\gamma$ is Euler's constant is part of this family. The integrals of
type (\ref{class-0})  are linear combinations of
\begin{equation}
J_{n}(t) := \ift e^{- t x} \left( \ln x \right)^{n} \, dx. 
\label{J-def}
\end{equation}
\noindent
The values of these integrals are expressed in terms of the gamma function
\begin{equation}
\Gamma(s) = \ift x^{s-1}e^{-x} \, dx 
\label{gamma-def}
\end{equation}
\noindent
and its derivatives.

\section{The evaluation} \label{evaluation}
\setcounter{equation}{0}

In this section we consider the value of $J_{n}(t)$ defined in (\ref{J-def}). 
The change of variables $s = tx$ yields
\begin{equation}
J_{n}(t) = \frac{1}{t} \ift e^{-s} \left( \ln s - \ln t \right)^{n} \, ds. 
\end{equation}
\noindent
Expanding the power yields $J_{n}$ as a linear combination of 
\begin{equation}
I_{m} := \ift e^{-x} \left( \ln x \right)^{m}\, dx, \quad 0 \leq m \leq n. 
\end{equation}

An analytic expression for these integrals can be obtained directly from 
the representation of the {\em gamma function} in (\ref{gamma-def}).

\begin{Prop}
For $n \in \mathbb{N}$ we have
\begin{equation}
\ift \left( \ln x \right)^{n} \, x^{s-1} e^{-x} \, dx = 
\left( \frac{d}{ds} \right)^{n} \Gamma(s). 
\end{equation}
\noindent
In particular
\begin{equation}
I_{n}:= \ift \left( \ln x \right)^{n} \,  e^{-x} \, dx = \Gamma^{(n)}(1). 
\label{iofn}
\end{equation}
\end{Prop}
\begin{proof}
Differentiate (\ref{gamma-def}) $n$-times with respect to the parameter 
$s$. 
\end{proof}

\medskip

\begin{Example}
Formula $\mathbf{4.331.1}$ in \cite{gr} states 
that\footnote{The table uses $C$ for 
the Euler constant.}
\begin{equation}
\ift e^{- \mu x} \, \ln x \, dx = - \frac{\delta}{\mu}
\end{equation}
\noindent
where $\delta = \gamma + \ln \mu$. This 
value follows directly by the change of  variables $s = \mu x$ and the 
classical special value $\Gamma'(1) = - \gamma$.  The reader will find in 
chapter 9 of \cite{irrbook} details on this constant.  In 
particular, if $\mu =1$, then
$\delta = \gamma$ and we obtain (\ref{class1}):
\begin{equation}
\ift e^{-x} \, \ln x \, dx = - \gamma. 
\end{equation}
\noindent
The change of variables $x = e^{-t}$ yields the form
\begin{equation}
\ifft t \, e^{-t} \, e^{-e^{-t}} \, dt = \gamma. 
\end{equation}
\end{Example}

\medskip

Many of the evaluations are given in terms of the {\em polygamma function}
\begin{equation}
\psi(x) = \frac{d}{dx} \ln\Gamma(x). 
\end{equation}
\noindent
Properties of $\psi$ are summarized in Chapter 1 of \cite{srichoi}. 
A simple representation is 
\begin{equation}
\psi(x) = \lim\limits_{n \to \infty} \left( \ln n - 
\sum_{k=0}^{n} \frac{1}{x+k} \right), 
\label{psi-def}
\end{equation}
\noindent
from where we conclude that
\begin{equation}
\psi(1) = \lim\limits_{n \to \infty} \left( \ln n - 
\sum_{k=1}^{n} \frac{1}{k} \right) = - \gamma, 
\end{equation}
\noindent
this being the most common definition of the Euler's constant $\gamma$. This 
is precisely the identity $\Gamma'(1) = - \gamma$. 

The derivatives of $\psi$ satisfy
\begin{equation}
\psi^{(m)}(x) = (-1)^{m+1} m! \, \zeta(m+1,x), \label{der-poly}
\end{equation}
\noindent
where 
\begin{equation}
\zeta(z,q) := \sum_{n=0}^{\infty} \frac{1}{(n+q)^{z}}
\end{equation}
\noindent
is the {\em Hurwitz zeta function}. This function appeared in \cite{moll-gr1} 
in the evaluation of some logarithmic integrals. 

\medskip

\begin{Example}
Formula $\mathbf{4.335.1}$ in \cite{gr} states that 
\begin{equation}
\ift e^{- \mu x} \left( \ln x \right)^{2} \, dx = 
\frac{1}{\mu} \left[ \frac{\pi^{2}}{6} + \delta^{2} 
\right], 
\label{form2}
\end{equation}
\noindent
where $\delta = \gamma + \ln \mu$ as before. This 
can be verified using the procedure described above: the change of 
variable $s = \mu x$ yields
\begin{equation}
\ift e^{- \mu x} \left( \ln x \right)^{2} \, dx = 
\frac{1}{\mu} \left( I_{2} - 2 I_{1} \ln \mu  + 
I_{0} \ln^{2}\mu \right), 
\end{equation}
\noindent 
where $I_{n}$ is defined in (\ref{iofn}). To complete the evaluation we 
need some special values: $\Gamma(1) = 1$ is elementary,  
$\Gamma'(1) = \psi(1) = - \gamma$ appeared above and using 
(\ref{der-poly}) we have
\begin{equation}
\psi'(x) = \frac{\Gamma''(x)}{\Gamma(x)} - 
\left( \frac{\Gamma'(x)}{\Gamma(x)} \right)^{2}. 
\end{equation}
\noindent
The value 
\begin{equation}
\psi'(1) = \zeta(2) = \frac{\pi^{2}}{6}, 
\end{equation}
\noindent
where $\zeta(z) = \zeta(z,1)$ is the Riemann zeta function, comes
directly from (\ref{der-poly}). Thus 
\begin{equation}
\Gamma''(1) = \zeta(2) + \gamma^{2}.
\end{equation}

Let $\mu=1$ in (\ref{form2}) to produce
\begin{equation}
\ift e^{-x} \left( \ln x \right)^{2} \, dx = \zeta(2) + \gamma^{2}.
\end{equation}

\medskip

Similar arguments yields formula $\mathbf{4.335.3}$ in \cite{gr}:
\begin{equation}
\ift e^{-\mu x} \left( \ln x \right)^{3} \, dx = - \frac{1}{\mu} 
\left[ \delta^{3} + \tfrac{1}{2} \pi^{2} \delta - \psi''(1) \right], 
\end{equation}
\noindent
where, as usual, $\delta = \gamma + \ln \mu$. The special 
case $\mu=1$ now yields
\begin{equation}
\ift e^{-x} \left( \ln x \right)^{3} \, dx = 
-\gamma^{3} - \tfrac{1}{2}\pi^{2}\gamma + \psi''(1).
\end{equation}
\noindent
Using the evaluation
\begin{equation}
\psi''(1) = -2 \zeta(3)
\end{equation}
\noindent
produces
\begin{equation}
\ift e^{-x} \left( \ln x \right)^{3} \, dx = 
-\gamma^{3} - \tfrac{1}{2}\pi^{2}\gamma -2 \zeta(3).
\end{equation}

\end{Example}

\medskip

\begin{Problem}
In \cite{irrbook}, page 203, we introduced the notion of {\em weight} for 
some real numbers. In particular, we have assigned $\zeta(j)$ the weight $j$. 
Differentiation increases the weight by $1$, so that $\zeta'(3)$ has 
weight $4$. The task is to check that the integral 
\begin{equation}
I_{n} := \ift e^{-x} \left( \ln x \right)^{n} \, dx 
\end{equation}
\noindent
is a homogeneous form of weight $n$. 
\end{Problem}

\section{A small variation} \label{small-variation}
\setcounter{equation}{0}

Similar arguments are now employed to produce a larger family of integrals. 
The representation
\begin{equation}
\ift x^{s-1} e^{- \mu x} \, dx = \mu^{-s} \Gamma(s),
\end{equation}
is differentiated $n$ times with respect to the parameter $s$ to produce
\begin{equation}
\ift \left( \ln x \right)^{n} x^{s-1} e^{-\mu x} \, dx = 
\left( \frac{d}{ds} \right)^{n} \left[ \mu^{-s} \Gamma(s) \right]. 
\end{equation}
\noindent
The special case $n=1$ yields
\begin{eqnarray}
\ift x^{s-1} e^{-\mu x} \, \ln x \, dx & =  &
\frac{d}{ds} \left[ \mu^{-s} \Gamma(s) \right] \label{special}  \\
& = & \mu^{-s} \left( \Gamma'(s) - \ln \mu \, \Gamma(s) \right) \nonumber \\
& = & \mu^{-s} \Gamma(s) \left( \psi(s) - \ln \mu \right). \nonumber
\end{eqnarray}
\noindent
This evaluation appears as $\mathbf{4.352.1}$ in 
\cite{gr}. The special case $\mu=1$
yields 
\begin{equation}
\ift x^{s-1} e^{-x} \, \ln x \, dx = \Gamma'(s), 
\end{equation}
\noindent
that is $\mathbf{4.352.4}$ in \cite{gr}. 

Special values of the gamma function and its derivatives yield more 
concrete evaluations. For example, the functional equation 
\begin{equation}
\psi(x+1) = \psi(x) + \frac{1}{x},
\end{equation}
\noindent
that is a direct consequence of $\Gamma(x+1) = x \Gamma(x),$ yields
\begin{equation}
\psi(n+1) = - \gamma + \sum_{k=1}^{n} \frac{1}{k}. 
\end{equation}
\noindent
Replacing $s=n+1$ in (\ref{special}) we obtain 
\begin{equation}
\ift x^{n} e^{-\mu x} \, \ln  x \, dx = 
\frac{n!}{\mu^{n+1}} \left( \sum_{k=1}^{n} \frac{1}{k} - \gamma - \ln \mu 
\right),
\end{equation}
\noindent
that is $\mathbf{4.352.2}$ in \cite{gr}. 

The final formula of Section $\mathbf{4.352}$ in \cite{gr} is $\mathbf{4.352.3}$
\begin{equation}
\ift x^{n-1/2} e^{- \mu x} \, \ln x \, dx = 
\frac{\sqrt{\pi} \, (2n-1)!!}{2^{n} \, \mu^{n+1/2}} 
\left[ 2 \sum_{k=1}^{n} \frac{1}{2k-1} - \gamma - \ln(4 \mu) \right]. 
\nonumber
\end{equation}
\noindent
This can also be obtained from (\ref{special}) by using the classical 
values
\begin{eqnarray}
\Gamma( n + \tfrac{1}{2} ) & = & \frac{\sqrt{\pi}}{2^{n}} (2n-1)!! \nonumber \\
\psi(n+ \tfrac{1}{2}) & = & -\gamma + 2 \left( \sum_{k=1}^{n} \frac{1}{2k-1} 
- \ln 2 \right). \nonumber
\end{eqnarray}
\noindent
The details are left to the reader. \\

Section $\mathbf{4.353}$ 
of \cite{gr} contains three peculiar combinations of integrands.
The first two of them can be verified by the 
methods described above: formula $\mathbf{4.353.1}$
states 

\begin{equation}
\ift (x - \nu) x^{\nu -1 } e^{-x} \, \ln x \, dx = \Gamma(\nu), 
\end{equation}

\noindent
and $\mathbf{4.353.2}$ is 

\begin{equation}
\ift (\mu x - n - \tfrac{1}{2}) x^{n- \tfrac{1}{2}} e^{-\mu x} \, 
\ln x \, dx = \frac{(2n-1)!!}{(2 \mu)^{n}} \sqrt{\frac{\pi}{\mu}}. 
\end{equation}
\noindent

\bigskip

\noindent
{\bf Acknowledgments}. The author wishes to thank Luis Medina for a  
careful reading of an earlier version of the paper. The partial support of 
$\text{NSF-DMS } 0409968$ is also acknowledged. 

\bigskip


\end{document}